\documentclass[preprint,12pt]{elsarticle}

\usepackage{amssymb}
\usepackage{amsmath}
\usepackage{amsthm}
\usepackage{hyperref}

\usepackage{siunitx}

\journal{Computers and Chemical Engineering}

\usepackage{xcolor}

\graphicspath{{./figures/}}

\begin{document}

\begin{frontmatter}



\title{Tracking in-silico Lagrangian sensors in a lab-scale stirred tank reactor}
\tnotetext[mytitlenote]{This project is funded by the Deutsche Forschungsgemeinschaft (DFG, German Research Foundation) – SFB 1615 – 503850735.}

\author[label1]{Vamika Rathi\corref{cor}}  
\author[label1]{Fatima Sehar}
\author[label1]{Finn Sommer}
\author[label1]{Sebastian Götschel}
\author[label2]{Eike Steuwe}
\author[label2]{Alexandra von Kameke}
\author[label1]{Daniel Ruprecht}
\affiliation[label1]{organization={Chair Computational Mathematics, Institute of Mathematics, Hamburg University of Technology},
             city={Hamburg},
             country={Germany}}
\affiliation[label2]{organization={Heinrich-Blasius Institute, Faculty of Engineering and Computer Science, Hamburg University of Applied Sciences},
             city={Hamburg},
             country={Germany}}

\cortext[cor]{Corresponding author}


\begin{abstract}
Lagrangian sensors have shown promise to improve operator awareness of conditions inside a chemical reactor but three-dimensional tracking remains a mostly unsolved challenge.
We explore a setup where in-silico sensors, based on a recently proposed real-world design, are tracked using data from an accelerometer and magnetometer available from a built-in inertial measurement unit.
Filtering algorithms, using a bespoke dynamical model, are used to process these readings into position estimates.
We compare tracking performance of an extended Kalman filter, a particle filter and the unscented Kalman filter implemented in the \textsc{pykalman} library. 
Our numerical experiments track in-silico particles moving in an analytically given three-dimensional vortex as well as in the experimentally measured flow-field of a lab-scale stirred tank reactor.
Using the Maxey-Riley-Gatignol equations for the movement of inertial particles as ground-truth, we demonstrate that  trajectories can be reconstructed from noisy synthetic data with errors below 10 \%.
\end{abstract}

\begin{graphicalabstract}
\end{graphicalabstract}

\begin{highlights}
\item We propose a framework to track Lagrangian sensors based on accelerometer and magnetometer readings
\item An extended Kalman filter and particle filter are proposed to fuse sensor readings with a bespoke dynamic model
\item Experiments for in-silico replicas of a recent Lagrangian sensor design show that trajectories can be reconstructed with relative errors below 5\%
\end{highlights}


\end{frontmatter}



\section{Introduction}\label{sec:intro}
Lagrangian sensors (LS) are sensors that travel freely in a liquid and provide measurements along their trajectory.
By contrast, Eulerian sensors are fixed and provide measurements over time at a single point in space.
LS are used in a variety of applications.
Examples are drifters in oceanography~\cite{AmadorEtAl2016}, coastal engineering~\cite{FeddersenEtAl2024} or to design pipes that avoid injury to fishes~\cite{CoxEtAl2026}.
In chemical process engineering, Lagrangian sensors have been explored for decades as means to provide data about the interior state of chemical reactors, starting with the concept of a ``radio pill'' proposed in 1969~\cite{Bryant1969,MannEtAl1981}.
Eulerian sensors, installed at the walls of a vessel, only give an incomplete picture of the reactions happening inside, in particular in bio-reactors where homogeneous mixing is difficult to achieve because too vigorous stirring would damage the cells~\cite{Bisgaard2021}.
Various devices have been proposed and studied over the years, see Bisgaard et al.~(2020) for a comprehensive survey and their Figure~2 for a timeline of developments~\cite{BisgaardEtAl2020}.
LS have also been conceptualized as an important part of a chemical reactor's digital twin~\cite{MuldbakEtAl2022}.
A broader perspective on flow-following sensor technology has emphasized its potential for validating CFD simulations and characterizing industrial flows~\cite{HaringaEtAl2025}.

LS have been equipped with a range of sensors to measure temperature, pressure, acceleration, pH, dissolved oxygen or magnetic fields~\cite{BisgaardEtAl2020}.
There are demonstrations in practical settings that show that LS can provide a range of useful information about the interior state of a chemical reactor, for example average residence times in vertical compartments, but collected data remain statistical rather than spatially or temporally resolved~\cite{Bisgaard2021,BisgaardEtAl2021,BisgaardEtAl2021b,ReineckeEtAl2022,HofmannEtAl2024,HofmannEtAl2025}.
A key obstacle is that localization remains a challenge.
In experimental setups, it is possible to track particles by some distinct property like color in video tracking, an emitted radio-signal, electrical resistance tomography, magnetic permeability or radioactivity in positron emission particle tracking (PEPT)~\cite{ReineckeEtAl2012,ReineckeEtAl2017,RautenbachEtAl2024,SykesEtAl2025}.
While extremely useful in a lab setting, these techniques cannot be used in standard operation.
Pressure readings and the assumption of hydro-static pressure can provide approximate vertical positioning and allow for the determination of averaged quantities with some accuracy~\cite{BisgaardEtAl2021,ReineckeEtAl2012,ReineckeEtAl2017} but do not provide information about horizontal position.

Fully time and space resolved tracking is complicated by the fact that commonly sized LS do not truly follow the flow but experience inertial effects depending on their size and velocity gradients in the fluid~\cite{HofmannEtAl2022}.
How closely a LS follows the flow is determined by its Stokes number, a non-dimensional parameter that represents that ratio of the particle response time to the characteristic time scale of the fluid~\cite{HofmannEtAl2022}.
Bisgaard et al.~(2022) show that LS can range in Stokes numbers from $\mathcal{O}(10^{-2})$ up to $\mathcal{O}(10)$, depending on their size and the viscosity of the carrier liquid.
There is evidence that particles with Stokes number $\mathcal{O}(10^{-1})$ and smaller can reasonably be considered flow-following with errors of $1 \%$ or less~\cite[\S 5.3.1]{Tropea2007}.
Because of the required instrumentation, however, particular LS with broader capabilities are larger and fall outside of this range.

Some sensor readings in an LP can be used to track its position, in particular readings from an IMU, magnetometer or pressure sensor.
Acceleration data from an inertial measurement unit IMU can be used for dead reckoning, but, on its own, suffers from rapidly accumulating drift errors over time.
Magnetic fields, either from the Earth or an artificial source, can act as beacons, and their measurements can be fused with IMU data for more robust tracking.
A magnetic field has been used, for example, to detect when a particle is at the apex of a chemical reactor to derive circulation times~\cite{ReineckeEtAl2017}.
Buntkiel et al.~(2023) study the performance of triaxial accelerometer in combination with a magnetometer  measuring Earth's magnetic field fused via a discrete-time Kalman filter for positioning~\cite{BuntkielEtAl2023}.
They validate reconstructed mean radial accelerations against measurements and, in their conclusions, propose to use an artificial time-dependent magnetic field to improve performance.
The idea of using an artificial magnetic field to track a LS called Sens-O-Sphere was proposed by Lange~\cite{Lange2019} but seems to have not yet been realized.

Motion tracking using readings from magnetometers and inertial sensors has been proposed to reduce computational requirements while maintaining orientation accuracy through EKF-based correction steps~\cite{CoxEtAl2023}.
An IMU/magnetometer-based indoor positioning system demonstrated that inertial data can be fused with magnetic field information to reduce drift in portable tracking systems~\cite{HellmersEtAl2013}.
Broader reviews of motion capture technologies in industrial environments highlight a growing interest in tracking freely moving bodies, but note that most systems rely on external references and kinematic motion models rather than physics-based dynamics~\cite{MeolottoEtAl2020}.
Wu et al. showed that IMU data fused with readings from an artificial dipole magnetic field via an extended Kalman filter can provide robust indoor positioning of a mobile cart~\cite{wu2016robust}.

In this paper, we follow up on the suggestion by Buntkiel et al.~~\cite{BuntkielEtAl2023} to use an artificial magnetic field and demonstrate that such a setup can track in-silico LP in the experimentally determined flow field of a lab-scale stirred tank reactor using a simulated artificial magnetic field.
In contrast to previous approaches, we use a physics-based state transition model in an extended Kalman filter (EKF) to account for the added mass and drag that LPs experience.
For validation, we employ the same dynamical model in a particle filter and compare its performance against that of the EKF.
To illustrate the benefit from using a physics-based transition model, we also compare tracking performance against the unscented Kalman filter (UKF) from the \textsc{pykalman}~\cite{pykalman} library.
All filters use the same dynamics and measurement models.
All the tracking results shown in the paper can be reproduced using the published code~\cite{Rathi_tracking}.

\section{Methodology}\label{sec:method}
In \S\ref{subsec:tracking} we describe the proposed tracking setup, including the physical characteristics of our in-silico particles.
\S\ref{subsec:marge} introduces the equations governing inertial particle motion that we use to create ground-truth trajectories.
In reduced form, they also serve as dynamic models used for tracking.
\S\ref{subsec:flows} describes the two velocity field required as input to the dynamical equations described in \S\ref{subsec:marge}, a 3D vortex and an experimental flow field of a lab-scale stirred tank reactor. 
The experimental setup that was used to obtain the second velocity field is also briefly described.
Finally, \S\ref{subsec:kalman} introduces the used extended Kalman filter while~\S\ref{subsec:particle} introduces the particle filter, including the necessary modifications to ensure accurate tracking like adaptive tempering and innovations. 
\subsection{Tracking framework}\label{subsec:tracking}
\begin{figure}[t]
	\centering
	\includegraphics{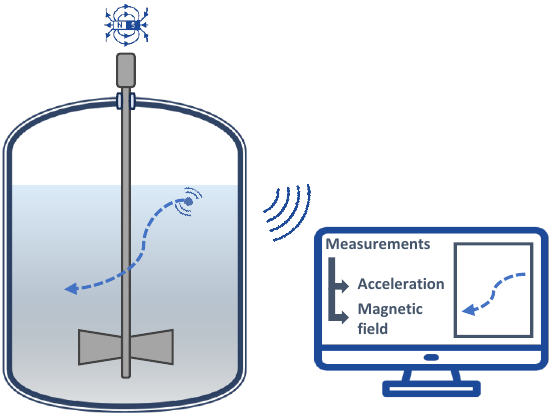}
	\caption{Illustration of a sensor inside a stirred tank reactor and position of the magnetic field source.}
	\label{fig:reactor}
\end{figure}

Figure~\ref{fig:reactor} shows a sketch of the proposed setup.
A magnetic field source placed above the reactor and generates a magnetic dipole field
\begin{equation} \label{mag_fld}
  B(\mathbf{r}) = \frac{\mu_0}{4\pi} \left( \frac{3 (\mathbf{m} \cdot \mathbf{r} ) \mathbf{r} - \mathbf{m} r^2}{r^5} \right)
\end{equation}
which is commonly used for magnetic field base positioning systems~\cite{PaskuEtAl2017}.
Here, $\mu_0$ is the magnetic permeability of the free space, $\mathbf{r}$ is the position vector originating from the magnetic field source, $r = \left\| \mathbf{r} \right\|$ the distance and $\mathbf{m}$ the magnetic moment of the coil generating the field, which is orthogonal to the coil surface.
Our in-silico particles are based on the real-world design proposed by Geläschus et al.~(2026)~\cite{GelaeschusEtAl2026}.
They contain an IMU with an accelerometer and vector magnetometer that can measure the strength and direction of the magnetic field~\cite{Edelstein2007}.
Data transmission will be realized using low-power communication methods suited to a reactor environment~\cite{GelaeschusEtAl2026} but, for the purpose of this paper, is assumed to work reliably except for artificially added Gaussian noise.

Note that many LS also feature pressure sensors that can be used for approximate vertical positioning by assuming hydrostatic pressure in the reactor vessel~\cite{BisgaardEtAl2021}.
Introducing this as an additional source of information into the filtering algorithms would be relatively straightforward, but requires some modifications and is thus left for future work.
While we test our approach for experimental flow-field data, our particles are in-silico, that is, simulated replicas, and the ground-truth trajectories we reconstruct are simulated.
However, their physical parameters summarized in Table~\ref{parameter-table} match those of the proposed design~\cite{GelaeschusEtAl2026}.
The corresponding fluid parameters of the experimental water-based setup are also included in Table~\ref{parameter-table}. 
\begin{table}[h]
	\centering
		\begin{tabular}{|l|c|c|}
			\hline
			Parameter & Symbol & Value \\
			\hline
			Particle radius & $a$ & \SI{0.0025}{\meter}\\
			Particle density & $\rho_p$ & \SI{1010}{\kilogram\per\meter^3}\\
			Fluid density & $\rho_f$ & \SI{998}{\kilogram\per\meter^3}\\
			Fluid kinematic viscosity & $\nu$ & \SI{1.004e-6}{\meter\squared\per\second}\\
			Fluid characteristic length & $L$ & \SI{0.130}{\meter}\\
			Fluid characteristic velocity & $U$ & \SI{0.07}{\meter\per\second}\\ 
			\hline	
			
		\end{tabular}
		\caption{Table of physical parameters of in-silico particles and the experimental flow field.}
	\label{parameter-table}
\end{table}

\subsection{Maxey-Riley-Gatignol Equations}\label{subsec:marge}
Because of their size, inertial effects are relevant for the motion of LS and they cannot be considered passive tracers~\cite{Bisgaard2021}.
Motion of spherical inertial particles in a fluid is described by the Maxey-Riley-Gatignol equation (MaRGE)~\cite{Maxey1983,Gatignol1983}.
Let $\boldsymbol{v} = \boldsymbol{\dot{x}}$ be the absolute velocity of a particle at position $\boldsymbol{x}$ with radius $a$ and density $\rho_p$.
MaRGE in dimensionless form then reads
\begin{equation}\label{MaRGE}
	\frac{\mathrm{d}\,\boldsymbol{v}}{\mathrm{d}t} =
	R\,\frac{\mathrm{D}\,\boldsymbol{u}}{\mathrm{D}t}
	- \frac{R}{S}(\boldsymbol{v}-\boldsymbol{u})
	- R\sqrt{\frac{3}{S\pi}}
	\frac{\mathrm{d}}{\mathrm{d}t}
	\int_{t_0}^{t} \frac{1}{\sqrt{t-\tau}}
	(\boldsymbol{v}-\boldsymbol{u})\,\mathrm{d}\tau
	- (1-R)\boldsymbol{G},
\end{equation}
using a characteristic time $T$ and characteristic velocity $U$ to non-dimensionalize, following the approach by Daitche~\cite{Daitche2013}.
Here, $t_0$ is the initial time, $\boldsymbol{u}(\boldsymbol{x}, t)$ is the fluid velocity at the particle position and $R$, $S$ and $\boldsymbol{G}$ are the dimensionless parameters
\begin{equation}
	R = \frac{3 \rho_f}{\rho_f + 2 \rho_p},  \; \; \; \; \; \; \; \; \; \; S = \frac{1}{3} \frac{a^2/\nu}{T},  \; \; \; \; \; \; \; \; \; \; \boldsymbol{G} = \frac{T}{U} \boldsymbol{g},
\end{equation}
where $\rho_f$ is the density of the fluid, $\nu$ is the kinematic viscosity of the fluid and $\boldsymbol{g}$ is the gravitational acceleration vector.
The material derivative
\begin{equation}
	\label{eq:mat-dev}
	\frac{\mathrm{D}\,\boldsymbol{u}}{\mathrm{D} t} = \frac{\partial\,\boldsymbol{u}}{\partial t} + \,\boldsymbol{u} \cdot \nabla \,\boldsymbol{u}
\end{equation}
in~\eqref{MaRGE} is calculated analytically for the vortex field.
For the experimental field, where only discrete data points are available, the gradient is computed by spline interpolation and the temporal derivative by linear interpolation.

Using the parameters in Table~\ref{parameter-table} and a characteristic time scale of $T = L / U \approx 1.86s$, we get $R = 0.992$ and $S = 1.12$.
Throughout the paper, we refer to the parameter $S$ as the Stokes number as it compares the particle relaxation time $a^2/\nu$ to the flow timescale $T$~\cite{UrizarnaEtAl2025} but note that slightly different uses exist in the literature.
We solve~\eqref{MaRGE} numerically using a linear multi-step method developed by Daitche~\cite{Daitche2013} and generalized to 3D by Rathi et al.~\cite{RathiEtAl2026}.
For various simplified flow fields, analytical solutions for $S$ ranging from 0.01 to 1 are provided by Prasath et al.~(2019)~\cite{PrasathEtAl2019}.

The purpose of this paper is to demonstrate that the proposed framework in principle allows tracking a particle moving inside the fluid.
We use~\eqref{MaRGE} including the Basset force modelled by the integral to generate ground-truth trajectories and synthetic data, to which artificial 5\% Gaussian noise is added.
The two filtering algorithms introduced below use a continuous process model that ignores the Basset force to reconstruct the ground-truth trajectories from the noisy synthetic data.
In essence, our setup introduces a controlled model error for in-silico particles.
By contrast, validation for real-world particles, which is planned for future work, would face additional challenges as the model error when using MaRGE to predict inertial particle motion is not yet fully understood which would make the physics-based dynamic models in the filters less accurate.

\subsection{Velocity fields}\label{subsec:flows}
The key assumption in MaRGE is that the impact of the particle's movement on the carrier liquid is confined to the integral term in~\eqref{MaRGE}, which models wake effects and accounts for temporal delays in boundary layer development~\cite{CroweEtAl1998}.
Therefore, the fluid flow field $\mathbf{u}$ in~\eqref{MaRGE} needs to be prescribed and is not modelled.
We consider two different flow fields, an analytically given 3D vortex and an experimentally determined velocity field in a lab-scale stirred-tank reactor.

\paragraph{Vortex flow}
Our 3D vortex is an extension of a simple two-dimensional vortex proposed by Candelier~\cite{CandelierEtAl2004}
\begin{equation}\label{analy_fld}
	\boldsymbol{u}(x, y, z, t) = \omega(z, t)
	                             \begin{bmatrix}
	                             	 -y \\
	                             	 x \\
	                             	 0
	                             \end{bmatrix},
\end{equation}
with an angular frequency $\omega(z, t) = \omega_0 + \alpha \, \sin^2(z) \, \cos^2(t)$ that depends on $z$ and $t$. 
We use $\omega_0 = 4$ and $\alpha = 0.2$.
The $z$ component for this field is zero, which means there is no velocity in the vertical direction but the particle will still sink or rise due to gravity.

\paragraph{Lab-scale stirred tank reactor}
The experimental flow field is obtained from a lab scale reactor with diameter \SI{130}{\milli\meter} and height \SI{230}{\milli\meter}.
The data for the experimental flow field were provided by the group of Alexandra von Kameke, using a 3 litre stirred tank reactor with two Rushton turbines~\cite{WEILAND2023100448}. 
A \textit{shake-the-box}~\cite{Schanz2016} method was used to compute a Lagrangian representation of the flow field from a large number of images of passive tracer particles taken by high-speed cameras.
Velocities were mapped to an Eulerian grid using the scattered interpolant method implemented in Matlab~\cite{KamekeEtAl2026}. 
Data covers a period of 2 seconds with a temporal resolution of 0.002 seconds for a total of 1001 snapshots. 
For efficiency, we consider only every fifth snapshot and discard the rest.
Linear interpolation is used to produce velocity data between snapshots.
To obtain a continuous representation $\hat{\mathbf{u}}$ of the measured fluid field that can be evaluated for any particle position $(x,y,z)$, the data on the Eulerian grid is interpolated using cubic spline interpolation provided by the \textit{RegularGridInterpolator} of the SciPy interpolation sub-package.\footnote{\url{https://docs.scipy.org/doc/scipy/reference/interpolate.html}}
The spatial gradient in~\eqref{eq:mat-dev} is computed by using finite differences provided by the Python package \textit{findiff}~\cite{findiff}.

To obtain the interpolated fluid velocity at some arbitrary time $t$,  we use the interpolation objects $\hat{\mathbf{u}}(x,y,z,t_i)$ and $\hat{\mathbf{u}}(x,y,z,t_{i+1})$ provided by SciPy for two nearest times $t_{i} \leq t \leq t_{i+1}$ for which data is available.
The velocity $\hat{\mathbf{u}}(x,y,z,t)$ is then calculated by linear interpolation.
\begin{equation*}
    \hat{\mathbf{u}}(x,y,z,t) = \hat{\mathbf{u}}(x,y,z,t_i) + \frac{\hat{\mathbf{u}}(x,y,z,t_{i+1}) - \hat{\mathbf{u}}(x,y,z,t_i)}{t_{i+1} - t_{i}}(t - t_i).
\end{equation*}
The temporal partial derivative in~\eqref{eq:mat-dev} is approximated by 
\begin{equation*}
    \frac{\partial\boldsymbol{u}}{\partial t} \approx \frac{\hat{\mathbf{u}}(x,y,z,t_{i+1}) - \hat{\mathbf{u}}(x,y,z,t_i)}{t_{i+1} - t_{i}}.
\end{equation*}

\begin{figure}[t]
	\centering
	\begin{minipage}{0.5\textwidth}
		\centering
		\includegraphics[scale=1]{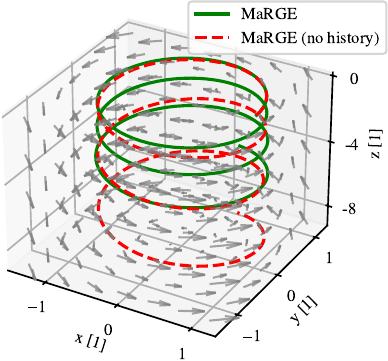}
		\label{Analytical}
	\end{minipage}%
	\begin{minipage}{0.5\textwidth}
		\centering
		\includegraphics[scale=1]{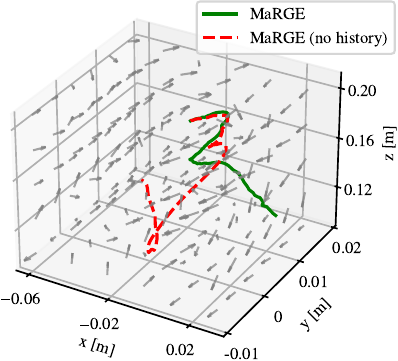}
		\label{Experimental}
	\end{minipage}
	\caption{The effect of history term on the trajectory of inertial particle.
		     The figure on the left is for the particle in analytical flow field for time $t \in [0,5]$.
		     The right figure is for the particle in experimental flow field form a stirred tank reactor for time $t \in [0,2] \sec$.
		     Trajectories obtained by MaRGE with (green, solid) and without (red, dashed) history term are plotted.}
	\label{fig:History-effect}
\end{figure}

\paragraph{Trajectories}
Figure~\ref{fig:History-effect} shows two trajectories for our in-silico 5$mm$ diameter device with the physical parameters given in Table~\ref{parameter-table}.
The device is denser than the fluid, which causes it to sink into the vortex, as shown in the left panel of Figure~\ref{fig:History-effect}.
In line with previous results, neglecting the Basset force modelled by the history term causes the particle to sink much faster.
For the experimental flow field, the impact of ignoring the Basset term is even more drastic.
Due to the lack of inertia, the particle moves at a completely different path.
This shows the ignoring the Basset force in the dynamic models of the filters introduces a controlled but significant model error which makes the reconstruction of the ground-truth from noisy synthetic data very challenging.

\subsection{Extended Kalman filter}\label{subsec:kalman}
This section describes how to adopt an Extended Kalman filter (EKF)\cite{mcgee1985discovery,becker2024kalman} for the described tracking setup.
EKF is the most commonly used nonlinear version of the Kalman Filter.
The nonlinear system equations are approximated through linearization about the current state estimate.
As process model we use $d \boldsymbol{x}(t)/dt = \boldsymbol{v}(t)$ and~\eqref{MaRGE} without the history term
\begin{equation}\label{MaRGEwoh}
	\frac{\mathrm{d}\,\boldsymbol{v}}{\mathrm{d}t} =
	R\,\frac{\mathrm{D}\,\boldsymbol{u}}{\mathrm{D}t}
	- \frac{R}{S}(\boldsymbol{v}-\boldsymbol{u})
	- (1-R)\boldsymbol{G}.
\end{equation}
Particle acceleration is denoted as $\boldsymbol{a} = \frac{\mathrm{d} \boldsymbol{v}}{\mathrm{d}t}$.
We write the system state at some time point $ t_k$ as $ \boldsymbol{S}_k = [\boldsymbol{x}_k \; \boldsymbol{v}_k \; \boldsymbol{a}_k]^T$, where $\boldsymbol{x}$, $\boldsymbol{v}$ and $\boldsymbol{a}$ are the 3D position, velocity and acceleration vectors of the LS.
Since explicit Euler is unconditionally unstable for second-order problems like~\eqref{MaRGEwoh}, we discretise the process model using symplectic Euler
\begin{subequations}
\label{eq:symplectic_euler}
\begin{align}
	\boldsymbol{v}_{k+1} & = \boldsymbol{v}_k + \Delta t \, \boldsymbol{a}_{k}, \\
	\boldsymbol{x}_{k+1} = \boldsymbol{x}_k + \Delta t \, \boldsymbol{v}_{k+1} & = \boldsymbol{x}_k + \Delta t \, \boldsymbol{v}_k + \Delta t^2 \, \boldsymbol{a}_{k},
\end{align}
\end{subequations}
where $\Delta t$ is the time step size and acceleration is computed from
\begin{equation}\label{MaRGEwohts}
		\boldsymbol{a}_{k} = R \frac{\mathrm{D}\,\boldsymbol{u}}{\mathrm{D} t} \bigg|_{t_k} - \frac{R}{S} (\boldsymbol{v}_k - \boldsymbol{u}_k) - (1-R)\boldsymbol{G}.
\end{equation}
This implementation allows to use a numerical time-step $\Delta t$ that is smaller than the sensor frequency, which allows for better control of the numerical error.
We use the following notations in the description of the EKF:
\begin{itemize}
	\centering
	\itemsep0em
	\item[$\boldsymbol{S}_{k+1,k}$] - predicted estimate of the system state at time $t_{k+1}$ calculated at time $t_k$.
	\item[$\boldsymbol{P}_{k+1,k}$] - predicted covariance of the system state at time $t_{k+1}$ calculated at time $t_k$.
	\item[$\boldsymbol{S}_{k+1,k+1}$] - updated estimate of the system state at time $t_{k+1}$.
	\item[$\boldsymbol{P}_{k+1,k+1}$] - updated covariance of the system state at time $t_{k+1}$.
\end{itemize}

We can now write the state extrapolation equation as
\begin{equation}
	\boldsymbol{S}_{k+1, k} = \begin{bmatrix}
								\boldsymbol{x}_{k+1,k} \\
								\boldsymbol{v}_{k+1,k} \\
								\boldsymbol{a}_{k+1,k}
							  \end{bmatrix}
								=
							\begin{bmatrix}
		                      \boldsymbol{x}_{k,k} + \Delta t \, \boldsymbol{v}_{k,k} + \Delta t^2 \, \boldsymbol{a}_{k,k} \\
		                      \boldsymbol{v}_{k,k} + \Delta t \, \boldsymbol{a}_{k,k} \\
		                      \boldsymbol{a}_{k,k}
                           \end{bmatrix} =: \boldsymbol{f}(\boldsymbol{S}_{k,k}).
\end{equation}
While this model appears linear in the state variables, the nonlinearity enters through the acceleration term $\boldsymbol{a}_{k,k}$, which is obtained from~\eqref{MaRGEwohts}. Specifically, the fluid velocity $\boldsymbol{u}(\boldsymbol{x},t)$ and its derivaties are nonlinear functions of the particle position $\boldsymbol{x}$ making the state transition function $\boldsymbol{f}$ nonlinear.
The covariance extrapolation reads
\begin{equation}
	\boldsymbol{P}_{k+1,k} = \boldsymbol{F}_{k} \boldsymbol{P}_{k,k} \boldsymbol{F}^{T}_{k} + \boldsymbol{Q},
\end{equation}
where $\boldsymbol{F}_k = \frac{\partial \boldsymbol{f}}{\partial \boldsymbol{S}}$ is the Jacobian of the state transition matrix and $\boldsymbol{Q}$ is the process noise covariance matrix
\begin{equation}
	\boldsymbol{Q} = \begin{bmatrix}
                       \Delta t^4 \boldsymbol{I}_3 & \Delta t^3 \boldsymbol{I}_3 & \Delta t^2 \boldsymbol{I}_3\\
                       \Delta t^3 \boldsymbol{I}_3 & \Delta t^2 \boldsymbol{I}_3 & \Delta t \, \boldsymbol{I}_3\\
                       \Delta t^2 \boldsymbol{I}_3 & \Delta t \, \boldsymbol{I}_3 & \boldsymbol{I}_3\\
	\end{bmatrix}\sigma^2_a,
\end{equation}
where $\boldsymbol{I}_3$ is the $3 \times 3$ identity matrix and $\sigma_a^2$ is the acceleration variance of the dynamic model.
The process noise covariance matrix is obtained for the dynamic model using the procedure by Becker~\cite{becker2024kalman}.
The extrapolation equations predict the state variables and the covariance.
Since we measure acceleration and magnetic field, the observation operator is
\begin{equation}
	\boldsymbol{h}(\boldsymbol{S}) = \begin{bmatrix}
		                \boldsymbol{a} \\
		                \boldsymbol{B}(\boldsymbol{x})
	\end{bmatrix},
\end{equation}
where $\boldsymbol{B}$ is the magnetic field from~\eqref{mag_fld}.
Thus, the measurement equation is
\begin{equation}
	\boldsymbol{\mathcal{Z}}_k = \boldsymbol{h}(\boldsymbol{S}_{k}).
\end{equation}
The measurement covariance matrix $\boldsymbol{R}_k$ is assumed to be constant
\begin{equation}
	\boldsymbol{R}_k = \begin{bmatrix}
		               \sigma^2_{a_m} \boldsymbol{I}_3 & \boldsymbol{0}_3 \\
		               \boldsymbol{0}_3 & \sigma^2_{B_m} \boldsymbol{I}_3 
	\end{bmatrix},
\end{equation}
where $\sigma^2_{a_m}$ is the acceleration measurement variance, $\sigma^2_{B_m}$ is the magnetic measurement variance and $ \boldsymbol{0}_3$ the $3\times3$ zero matrix.
The Kalman gain is calculated as
\begin{equation}
	\boldsymbol{K}_{k+1} = \boldsymbol{P}_{k+1,k} \, \frac{\partial \boldsymbol{h}}{\partial \boldsymbol{S}}^T \left(\frac{\partial \boldsymbol{h}}{\partial \boldsymbol{S}} \boldsymbol{P}_{k+1,k} \, \frac{\partial \boldsymbol{h}}{\partial \boldsymbol{S}}^T + \boldsymbol{R}_k \right)^{-1},
\end{equation}
where $\frac{\partial \boldsymbol{h}}{\partial \boldsymbol{S}}$ is the Jacobian of the nonlinear observation operator $\boldsymbol{h}(\boldsymbol{S})$ with respect to the system state.
Since the first three rows of $\boldsymbol{h}(\boldsymbol{S})$ are just acceleration components and the magnetic field only depends on positions, we have
\begin{equation}
	\frac{\partial \boldsymbol{h}}{\partial \boldsymbol{S}} = \begin{bmatrix}
		 \boldsymbol{0}_3 & \boldsymbol{0}_3 & \boldsymbol{I}_3 \\
	     \frac{\partial \boldsymbol{B}}{\partial \boldsymbol{x}} & \boldsymbol{0}_3 & \boldsymbol{0}_3	  
	\end{bmatrix}.
\end{equation}
The Jacobian $\frac{\partial \boldsymbol{B}}{\partial \boldsymbol{x}}$ of the dipole magnetic field model is computed via automatic differentiation (AD) using the forward mode AD implementation in JAX~\cite{jax2018github}.
Finally, the state update and covariance update equations, given the $(k+1)$th measurement $\boldsymbol{\tilde{\mathcal{Z}}}_{k+1}$, are
\begin{align}
	\boldsymbol{S}_{k+1,k+1} & = \boldsymbol{S}_{k+1,k} + \boldsymbol{K}_{k+1} (\boldsymbol{\tilde{\mathcal{Z}}}_{k+1} - \boldsymbol{h}(\boldsymbol{S}_{k+1,k})),\\
	\boldsymbol{P}_{k+1,k+1} & = \left( \boldsymbol{I}_9 - \boldsymbol{K}_{k+1} \frac{\partial \boldsymbol{h}}{\partial \boldsymbol{S}} \right) \boldsymbol{P}_{k+1,k} \left( \boldsymbol{I}_9 - \boldsymbol{K}_{k+1} \frac{\partial \boldsymbol{h}}{\partial \boldsymbol{S}} \right)^T + \boldsymbol{K}_{k+1} \boldsymbol{R}_{k+1} \boldsymbol{K}_{k+1}^T.
\end{align}

\subsection{Particle filter}\label{subsec:particle}
For comparison and cross-validation, we compare the EKF against a particle filter.
Particle filters (PF) are sequential Monte Carlo methods for state estimation  that approximate the filtering distribution by associating weights with random samples (called particles, not to be confused with our in-silico LS)~\cite{DjuricPF2003, articleBF}. 
We implement a bootstrap PF which samples from the process model and uses the measurement likelihood as weights~\cite{Elfring2021ParticleFA}.
The posterior distribution  at time step $k$, modelling the probability of the particle being in a certain state, is represented by a set of $N_p$
weighted particles.
In contrast to the EKF, the PF does not include acceleration in the state explicitly so that the state vector in the PF is defined as $\boldsymbol{X}_k = [\boldsymbol{x}_k \; \boldsymbol{v}_k]^T$.
Let $\boldsymbol{X}_k^{(i)}$ denote the $i$-th particle's state with weight (probability) $w_k^{(i)}$.
The final state estimate delivered by the PF is the weighted mean
\begin{equation}
	\label{eq:pf_state_estimate}
	\hat{\boldsymbol{X}}_k = \sum_{i=1}^{N_p} w_k^{(i)} \boldsymbol{X}_k^{(i)}
\end{equation} 
of the particle ensemble.
Weights and particle states are continuously updated as new measurements become available.
 
Particles for the filter are initialized by sampling from a Gaussian distribution centred at the initial estimate $(\boldsymbol{x}_0, \boldsymbol{v}_0)$ with standard deviation $\sigma_{\text{init}} = 0.05$. 
The state propagation equations used to update the state of the particle with index $i$ are the same discretization of MaRGE
\begin{align}
\boldsymbol{v}_{k+1}^{(i)} &= \boldsymbol{v}_k^{(i)} + \Delta t \, \boldsymbol{a}_k^{(i)} + \boldsymbol{n}_{v,k}^{(i)}, \\
\boldsymbol{x}_{k+1}^{(i)} &= \boldsymbol{x}_k^{(i)} + \Delta t \, \boldsymbol{v}_{k+1}^{(i)} + \boldsymbol{n}_{x,k}^{(i)}, 
\end{align}
that is used in the EKF~\eqref{eq:symplectic_euler}, except for the addition of process noise $\boldsymbol{n}_{v,k}^{(i)}, \boldsymbol{n}_{x,k}^{(i)}  \sim \mathcal{N}(\boldsymbol{0}, \boldsymbol{Q})$ sampled independently for each particle $i$ at each time step $k$.
As in the EKF, acceleration is computed from the particle state by~\eqref{MaRGEwohts} and  the numerical timestep $\Delta t$  can be chosen to be smaller than the sensor frequency to reduce the discretization error.
The process noise covariance matrix 
\begin{align}
\boldsymbol{Q} = \begin{bmatrix}
\Delta t^4 \boldsymbol{I}_3 & \Delta t^3 \boldsymbol{I}_3 \\
\Delta t^3 \boldsymbol{I}_3 & \Delta t^2 \boldsymbol{I}_3
\end{bmatrix} \sigma_a^2.
\end{align}
is based on the discrete white-noise acceleration model by Becker~\cite{becker2024kalman}.

The PF maps the particle state to measurements via
\begin{align}
\boldsymbol{a}_{\text{pred}}^{(i)} &= \boldsymbol{a}_k^{(i)}, \\
\boldsymbol{B}_{\text{pred}}^{(i)} &= \boldsymbol{B}(\boldsymbol{x}_k^{(i)}).
\end{align}
where $\boldsymbol{a}_k$ is the acceleration as in~\eqref{MaRGEwohts} and $\boldsymbol{B}(\boldsymbol{x}_k^{(i)})$ is the magnetic field from~\eqref{mag_fld} at the particle position $\boldsymbol{x}_k^{(i)}$.

Because acceleration and magnetic field data differ in physical units and magnitude, we normalize both measurement residuals (also called innovations) by their respective uncertainties to ensure that they make comparable contribution to a particle's weight~\cite{becker2024kalman}.
Acceleration innovations are normalized by combining the empirical particle spread of the measurement noise
\begin{equation}
	\label{eq:z_acc_i}
\mathcal{Z}_{\text{acc},j}^{(i)} = \frac{a_{\text{pred},j}^{(i)} - \mathcal{Z}_{\text{acc},j}}{\sqrt{\sigma_{\text{ensemble},j}^2 + \sigma_{\text{meas}}^2}},
\end{equation}
where $\sigma_{\text{ensemble},j}$ is the standard deviation of particle predictions in component $j \in \{x, y, z\}$. This adaptive scaling down-weights acceleration measurements when particle predictions are inconsistent (large $\sigma_{\text{ensemble}}$), reflecting increased model uncertainty.

To ensure that magnetic field measurements remain valuable even when particles are spatially dispersed, magnetic field innovations are only scaled relative to field magnitude
\begin{equation}
	\label{eq:z_mag_i}
\mathcal{Z}_{\text{mag},j}^{(i)} = \frac{\boldsymbol{B}_{\text{pred},j}^{(i)} - \mathcal{Z}_{\text{mag},j}}{\sigma_{\boldsymbol{B}} \cdot |\mathcal{Z}_{\text{mag},j}|},
\end{equation}
where $\sigma_{\boldsymbol{B}}$ represents the relative measurement uncertainty of the magnetometer.

We assume conditional independence between sensor modalities and combine their contributions through multiplicative likelihood functions~\cite{Siebler2022,Zhang2011}. 
Particle weights are updated via importance sampling
\begin{equation}
w_k^{(i)} \propto w_{k-1}^{(i)}  p(\boldsymbol{\mathcal{Z}}_{\text{acc}}|X_k^{(i)})  p(\boldsymbol{\mathcal{Z}}_{\text{mag}}|X_k^{(i)})
\end{equation}
which becomes additive after applying a logarithm. 
Note that $\boldsymbol{\mathcal{Z}}_{\text{acc}}, \boldsymbol{\mathcal{Z}}_{\text{mag}} \in \mathbb{R}^3$ are vectors with components given by~\eqref{eq:z_acc_i} and~\eqref{eq:z_mag_i}.

\paragraph{Resampling}
Particle filters can suffer from ensemble degeneracy, when all but a very small number of particles have negligible importance weights.
This represents essentially a collapse of particle diversity and leads to a substantial amount of superfluous computation, where particles keep getting updated that do no longer contribute meaningfully to the overall estimate~\eqref{eq:pf_state_estimate}.
Particle weight degeneracy can be measured using the effective sample size
\begin{equation}
	\text{ESS} = \frac{1}{\sum_{i=1}^{N_p} (w_k^{(i)})^2}.
\end{equation}
We use adaptive tempering of the likelihood function to help maintain particle diversity~\cite{HerbstEtAl2019}
\begin{equation}
w_k^{(i)} \propto w_{k-1}^{(i)} \exp\left(\frac{\ell^{(i)}}{\tau}\right),
\end{equation}
where $\tau \geq 1$ is the adaptive temperature parameter and is doubled iteratively until ESS $\geq 0.5 N_p \text{(number of particles)}$ or $\tau = 64$.
Here, $\ell^{(i)}$ is a fusion parameter introduced as
\begin{equation}
	\label{eq:fusion}
	\ell^{(i)} = (1-\gamma) \ell_{\text{acc}}^{(i)} + \gamma \ell_{\text{mag}}^{(i)},
\end{equation}
with $\gamma \in [0,1]$ to explore different weighting between accelerometer and magnetometer readings.
However, in experiments not documented here we found that equal weighting by using $\gamma = 0.5$ works well and no tuning is required.
The log-likelihoods in~\eqref{eq:fusion} are computed from
\begin{align}
	\ell_{\text{acc}}^{(i)} &= -\frac{1}{2}\sum_{j=1}^{3} (\mathcal{Z}_{\text{acc},j}^{(i)})^2, \\
	\ell_{\text{mag}}^{(i)} &= -\frac{1}{2}\sum_{j=1}^{3} (\mathcal{Z}_{\text{mag},j}^{(i)})^2.
\end{align}
Note that filter performance can be sensitive to these thresholds and may require tuning for different applications.

Once the particle diversity becomes too low and we start seeing significant weight degeneracy, indicated by $\text{ESS} < 0.5 N_p$, we resample.
Resampling means generating a new set of more diverse particles to replace the old degenerate ensemble.
Kuptametee et al.~(2022) summarized and compared various resampling techniques for particle filters~\cite{KuptameteeEtAl2022}.
We use systematic resampling, which draws a single initial uniform random sample and generates the remaining particles at fixed intervals over equally divided segments of $[0,1]$, resulting in improved computational efficiency and reduced variance.
After resampling, we apply roughening by adding Gaussian noise to the particle state~\cite{LiEtAl2014}
\begin{align}
\boldsymbol{x}_i &\leftarrow \boldsymbol{x}_i + \mathcal{N}(\mathbf{0}, c^2\sigma_x^2 \mathbf{I}), \\
\boldsymbol{v}_i &\leftarrow \boldsymbol{v}_i + \mathcal{N}(\mathbf{0}, c^2\sigma_v^2 \mathbf{I}),
\end{align}
to restore particle diversity and prevent sample impoverishment.
We use $c = 0.15 N_p^{-1/6}$, $\sigma_x =$ \SI{1.0}{\meter}, and $\sigma_v =$ \SI{0.5}{\meter\per\second} for all experiments without further tuning as these values produce robust and accurate results.
Table~\ref{tab:pf_params} summarizes all particle filter parameters used in this study.
\begin{table}[h!]
\centering
\begin{tabular}{|l|c|c|}
\hline
Parameter & Symbol & Value \\
\hline
Number of particles & $N_p$ & 500 \\
Initial particle spread & $\sigma_{\text{init}}$ & \SI{0.05}{\meter}, \SI{0.05}{\meter\per\second}\\
Fusion weight & $\gamma$ & 0.5 \\
ESS resampling threshold & $0.5 N_p$ & 250 \\
Maximum tempering factor & $\tau_{\text{max}}$ & 64 \\
Roughening position scale & $\sigma_x$ & \SI{1.0}{\meter} \\
Roughening velocity scale & $\sigma_v$ & \SI{0.5}{\meter\per\second} \\
Substep factor & $f$ & 1 \\
\hline
\end{tabular}
\caption{Parameters used in the particle filter. Note that the particles in the filter correspond to points in state space indicating hypotheses about the state of the in-silico particle, but are not particles in the same sense.}
\label{tab:pf_params}
\end{table}

\section{Numerical results}\label{sec:results}
This section presents the results of tracking in-silico particles in the two flow fields, obtained using the EKF, PF and UKF.
For UKF, we use the additive unscented Kalman filter implementation from the \textsc{pykalman} library with the same initial state covariance, process noise covariance, and measurement noise covariance as for the EKF described in \S\ref{subsec:kalman}.
However, since the UKF needs positive-definite covariance matrices to calculate Cholesky decomposition~\cite{JulierEtAl1997}, we regularize these matrices by a small positive shift on the diagonal of $10^{-6}$.

The choice of initial covariance matrix ($\boldsymbol{P}_{0,0}$) is based on the confidence we have in the initial state.
We use the value $\boldsymbol{P}_{0,0} = 0.1 \boldsymbol{I}_{9\times9}$ for both vortex and experimental flows, reflecting moderate uncertainty in the initial state estimate.
The measurement noise covariances ($\sigma^2_{a_m}$ and $\sigma^2_{B_m}$) provided to the filters  are intentionally set to values different from the true noise levels, reflecting a realistic scenario where the exact sensor noise characteristics are not perfectly known.
Finally, the value for the acceleration variance of the dynamic model $\sigma^2_a$ was tuned to minimize the tracking error across all three filters.
The values of the parameters used for the filters are summarized in Table~\ref{tab:filter_parms}.

\begin{table}
	\centering
	\begin{tabular}{|l|c|c|}
		\hline
		Parameter & Symbol & Value \\
		\hline
		Acceleration measurement noise covariance &  $\sigma^2_{a_m}$ & 0.04 \\
		Magnetic measurement noise covariance &  $\sigma^2_{B_m}$ & 0.04 \\
		Dynamic model acceleration variance &  $\sigma^2_a$ & 0.8 \\
		\hline
	\end{tabular}
	\caption{Parameters of the noise added to the synthetic data fed to the three filters.}
	\label{tab:filter_parms}
\end{table}

Ground-truth trajectories and synthetic data are obtained by solving~\eqref{MaRGE} with high accuracy with zero initial relative velocity.
A $5\%$ Gaussian noise is added to all measurements before they are given to the filters.
To assess the accuracy of the reconstructed trajectories we use the relative error
\begin{equation}
	\label{error}
	\mathrm{RelErr}(t_k)
	= \frac{\lVert \hat{\boldsymbol{x}}(t_k) - \boldsymbol{x}(t_k) \rVert_2}{L},
	\qquad
\end{equation}
where $\lVert \cdot \rVert_2$ denotes the Euclidean norm, $\hat{\boldsymbol{x}}(t_k)$ and $\boldsymbol{x}(t_k)$ are the estimated and true particle positions at time $t_k$, and $L$ is the total arc length of the true trajectory.

The magnetic field source is placed at $(0,0,0.3)$ for both the analytical vortex and experimental cases.
In the analytical case, all quantities are treated as dimensionless from the outset, so the magnetic field model is evaluated directly at the nondimensional particle positions.
However, for the experimental case, the particle position is first redimensionalized by the characteristic length $L$ before evaluating the dipole field model~\eqref{mag_fld}, and the resulting field is subsequently nondimensionalized using the Earth's magnetic field as characteristic field strength $B_c = 2.5\times10^{-5}$ T.
Therefore, the corresponding measurement Jacobian is rescaled by a factor of $L/B_c$.

\subsection{Tracking in-silico particles in the vortex}

\begin{figure}[t!]
	\centering
	\begin{minipage}{0.5\textwidth}
		\centering
		\includegraphics[scale=1]{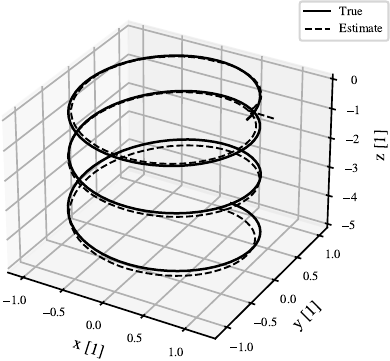}
	\end{minipage}%
	\begin{minipage}{0.5\textwidth}
		\centering
		\includegraphics[scale=1]{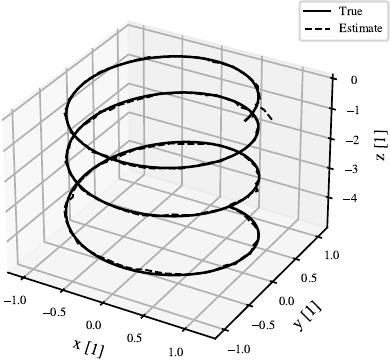}
	\end{minipage}
	\caption{Reconstructed trajectories over the time interval $[0,5]$ in the three-dimensional vortex. The solid line represents the ground-truth while the dashed line shows the estimates provided by the EKF (left) and the PF (right).}
	\label{fig:ana_results}
\end{figure}

Figure~\ref{fig:ana_results} shows the trajectories in the flow field~\eqref{analy_fld} reconstructed by the EKF and PF over the time interval $t \in [0, 5]$.
The true initial position of the particle is $(1,0,0)$ while the initial position state for the filters is  $(1.2,0.2,-0.1)$.
Both true and estimated initial relative velocity is set to zero.
For the EKF's and UKF, where the state includes acceleration, the initial acceleration is set to zero.

\begin{figure}[t!]
	\centering
	\includegraphics{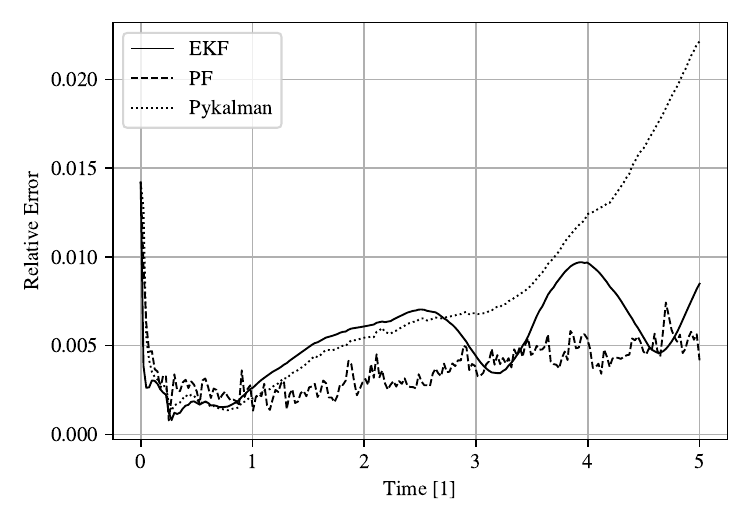}
	\caption{Relative error~\eqref{error} over time for the EKF, PF and UKF (Pykalman) for an in-silico LS in the three-dimensional vortex.}
	\label{fig:ana_comparison}
\end{figure}

Figure~\ref{fig:ana_comparison} shows the relative Euclidean error~\eqref{error} over time for the EKF, PF and UKF (Pykalman).
The EKF and PF achieve comparable accuracy and maintain a relative error below $1\%$.
By contrast, the UKF diverges progressively after $t \approx 2$ and produces errors more than three times larger by time $t = 5$.
All three filters show an initial spike around $t = 0$, because of the discrepancy between the true and assumed initial positions but can correct his as more and more sensor readings arrive.
It was also observed that the UKF is considerably more sensitive to the choice of the initial state covariance $\boldsymbol{P}_{0,0}$ than the EKF.


\subsection{Tracking in-silico particles in the lab-scale STR}
In the experiment, the velocity field was measured over a duration of 2 seconds and we track the in-silico particle over the full time interval.
Since MaRGE is solved numerically in nondimensional form, the experimental data must also be nondimensionalized.
This is done by scaling the velocity field with the characteristic velocity $U$ and the time domain $[0, 2]$ by the characteristic time $T = L/U$.
As a result, the trajectories obtained using MaRGE, as well as the EKF and PF estimates, are expressed in 
nondimensional form.
However, to allow comparison of the results shown here to the experimental setup, all shown quantities have been redimensionalized by multiplication with their respective reference.
The initial position of the particle is set to $(0, 0, 0.21)$~m while the initial position estimate for both filters differs from the true initial position by approximately $1\%$.

\begin{figure}[t]
	\centering
	\includegraphics{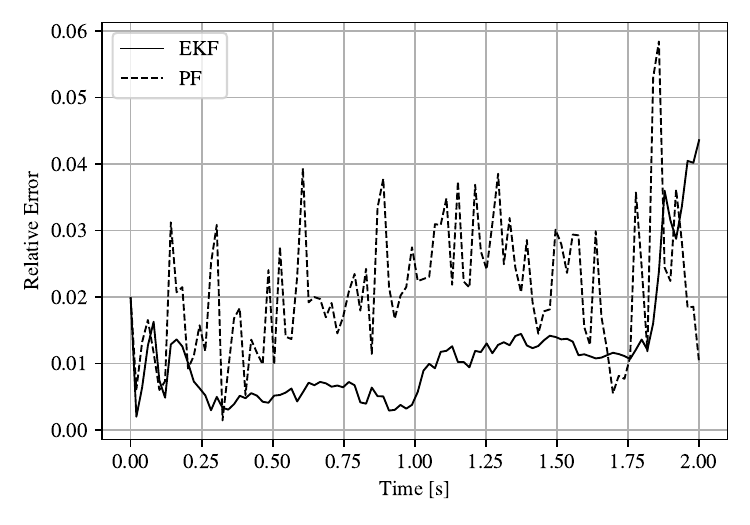}
	\caption{Relative Euclidean error~\eqref{error} over time for the extended Kalman filter and particle filter (PF) for an in-silico LS in the lab-scale STR.}
	\label{fig:exp_comparison}
\end{figure}

Figure~\ref{fig:exp_comparison} shows the relative Euclidean error over time for the EKF and PF in the experimental flow field.
The EKF produces smaller error than the PF throughout most of the tracking interval, although it produces a more accurate estimate at the very end.
Errors for both filters are below 6\% throughout the entire interval with averages around 1.1\% for the EKF and 2.2\% for the PF.
Due to the filter's stochastic nature, the PF error is considerably more volatile.
This suggests that the particle ensemble struggles to represent the posterior distribution accurately in the more complex, experimentally obtained flow field.
Both filters show a small initial transient near $t = 0$, consistent with the slight offset in the initial position estimate.
The UKF was also tested. 
However, it resulted in a mean relative error exceeding 1800\%, likely due to the strong nonlinearities present in the experimental flow field, and is therefore omitted from the comparison.
 

\section{Conclusions}\label{sec:conc}
The three-dimensional tracking of Lagrangian sensors in chemical reactors is challenging due to the complex dynamics.
Our paper proposes a tracking framework that combines filtering algorithms (EKF, UKF and PF) using physics-based models with sensor readings from an accelerometer and magnetometer measuring an artificial magnetic field.
The framework is tested on in-silico replicas of the sensors recently proposed by Geläschus et al.~\cite{GelaeschusEtAl2026}, first in an idealized three-dimensional vortex and then in a lab-scale stirred tank reactors, where the flow-field was characterized by particle image velocimetry.

We show that EKF and PF can reliably reconstruct trajectories of in-silico particles  with relative errors of the order of a few percent, even when faced with noisy data and an imperfect process model.
By contrast, the generic UKF provided by the pyKalman library works well for the vortex flow field but cannot cope with the nonlinearities present in the flow field of the stirred tank reactor.
While substantial challenges remain, see the discussion of limitations and next steps below, the proposed framework shows promise to eventually allow online three-dimensional tracking of Lagrangian sensors.

\paragraph{Limitations and future work}
While our study demonstrates the promise of combining IMU data with bespoke filters to enable full 3D tracking of Lagrangian sensors in chemical reactors, important challenges remain that need to be addressed before the approach could be used operationally.
Most importantly, while the flow field in the lab-scale reactor was measured experimentally, the reconstructed trajectories are simulated, even though neglecting the Basset force in the filters introduces a substantial modelling error that both EKF and PF were able to cope with.
The Maxey-Riley-Gatignol equations (MaRGE) are generally accepted as a useful model for the dynamics of inertial particles~\cite[\S 5.3.1]{Tropea2007}.
However, to date there is no analysis of how well MaRGE describes the movement of Lagrangian sensors, which are relatively large inertial particles and where the geometric center and the center of mass might not be perfectly aligned.
Therefore, we do not know the gap between the dynamic models used in the filters and the motion of real-world LS with high certainty.
Reconstructing experimentally measured trajectories of real Lagrangian particles in the lab-scale reactor is planned for future work.

Furthermore, our proposed setup, combining readings from an accelerometer and magnetometer for tracking using an externally generated magnetic field, can only be used in reactors made out of non-magnetic materials.
In steel vessels, for example, other sensors would have to be employed to counteract the drift from dead reckoning.
Pressure sensors has been used for positioning particles in the vertical~\cite{BisgaardEtAl2021} and could also be integrated into our filters for improved vertical positioning.
Other options might be tracking based on ultrasound~\cite{Chen2024}.
While our filter algorithms can be adopted to include additional or different sensor modalities, this would influence their tracking performance and thus require further tests and investigation.

In our tracking experiments, we assumed that sensor noise is additive Gaussian noise but the noise characteristics for real particles, having to transmit their IMU readings to a processing unit outside of the reactor via low-power communication, might be quite different.
However, modern machine learning architectures like autoencoders have proven to be quite effective at denoising sensor data and could be a useful remedy~\cite{Berahmand2024}.
Therefore, the issue of sensor and communication noise can likely be dealt with satisfactorily.

Lastly, the tracking algorithms in this study were run offline and cannot yet track particles in real-time.
Enabling online tracking will require to speed up computations by about two orders of magnitude.
This would require to translate the current Python code into a more efficient compiled programming language like C/C++ and to make effective use of compiler optimization as well as the use parallel processing on multi-core CPUs or graphics processing units (GPUs).
While this would require a substantial investment in terms of development time, the necessary techniques are well tested and readily available.

\bibliographystyle{elsarticle-num} 
\bibliography{refs}



\end{document}